\documentclass[12pt,letterpaper]{article}
\usepackage{url}

\usepackage{algpseudocode} 
\usepackage{algorithm}
\usepackage{amssymb}
\usepackage{graphicx}
\usepackage{fullpage}
\usepackage{amsfonts}
\usepackage{amsmath}
\usepackage{amsthm}
\usepackage{cite}

\begin{document}
\newtheorem{theorem}{Theorem}[section]
\newtheorem{lemma}[theorem]{Lemma}

\theoremstyle{definition}
\newtheorem{definition}[theorem]{Definition}
\newtheorem{example}[theorem]{Example}
\newtheorem{xca}[theorem]{Exercise}
\newtheorem{conjecture}[theorem]{Conjecture}

\theoremstyle{remark}
\newtheorem{remark}[theorem]{Remark}

\numberwithin{equation}{section}

\title{A short note on Jacobsthal's function}
\author{Fintan Costello and Paul Watts}

\date{}

\maketitle

\begin{abstract}
The function $g(n)$ represents the smallest number $Q$ such that every
sequence of $Q$ consecutive integers contains an integer coprime to
$n$.  We give a new and explicit upper bound on this function.
\end{abstract}

\section{Introduction}

Jacobsthal's function $g(n)$ represents the smallest number $Q$ such
that every sequence of $Q$ consecutive integers contains an integer
coprime to $n$.  This function has been studied by a number of
different authors, and is central to results on the maximal gaps
between consecutive primes \cite{MaierPomerance,Pintz} and on the
least prime in arithmetic progressions \cite{Pomerance}.  Taking $k$
to represent $\omega(n)$, the number of distinct prime factors of $n$,
Iwaniec's proof \cite{Iwaniec} that
\begin{equation*}
 g(n) \leq X \ (k \log  k )^{2}
\end{equation*} 
for some unknown constant $X$ gives the best asymptotic upper bound on
$g(n)$.  The best known explicit upper bounds, of
\begin{eqnarray*}
g(n) \leq 2^{k}&\mathrm{and}&g(n)\leq 2k^{2+2 e \log  k},
\end{eqnarray*}
are due to Kanold \cite{Kanold} and Stevens \cite{Stevens}
respectively, with the second bound being stronger for $k \geq 260$.
In this short note we describe  improvements on Stevens'
method that give a better bound of
$$g(n) \leq  2 e^{\gamma}  k^{5+5 \log \log k}$$
for $k>120$, where $\gamma$ is the Euler-Mascheroni constant.

We take $p$ to represent some prime and $p_i$ to represent the
$i{}^{\mathrm{th}}$ prime.  Because $g(n)=g(\mathrm{rad}(n))$, we may
assume that $n$ is squarefree.  For computational reasons, we assume
$k > 61$.  Define
$$P = \prod_{p \mid n} \left(1-\frac{1}{p}\right)$$
$$T_s =\sum_{r=1}^{s} (-1)^r \sum_{d \mid n, \omega(d)=r}\frac{1}{d}$$
$$T^{'}_{s} = T_s - P$$
where $s \geq 1$.  Given these definitions Stevens proves that if an odd value of $s$ is
chosen such that $P>T^{'}_{s}$, then $g(n) \leq Q$ when
\begin{equation}
\label{stevens}
 \frac{k^s}{P-T^{'}_{s}} <Q
\end{equation}
To determine suitable values of $s$, Stevens first notes that
\begin{equation*}
r! \sum_{d \mid n, \omega(d)=r}\frac{1}{d} <\left(\sum_{i=1}^{k}\frac{1}{p_i}\right)^r 
\end{equation*}
From \cite{RosserSchoenfeld} we have
\begin{equation*}
 \sum_{i=1}^{k}\frac{1}{p_i}  <  \log(C_k \log p_k)
\end{equation*}
for $k > 61$ where
$$C_k = e^{(B + \frac{1}{2 \log^2 p_k}) } $$
and  $B= 0.261498$, and so
\begin{equation*}
r! \sum_{d \mid n, \omega(d)=r}\frac{1}{d} < \left(\log(C_k \log p_k)\right)^r
\end{equation*}
for all $r \leq k$.  We thus have
\begin{equation*}
\begin{split}
T_s^{'} = &\sum_{r=s+1}^{k} (-1)^r \sum_{d \mid n,
  \omega(d)=r}\frac{1}{d} < \sum_{r=s+1}^{k} \frac{1}{r!}
\left(\sum_{i=1}^{k}\frac{1}{p_i}\right)^r \\ & <
\sum_{r=s+1}^{\infty} \frac{\left(\log(C_k \log p_k)\right)^r}{r!} <
\left(C_k \log p_k \right) \frac{ \left(\log(C_k\log
  p_k)\right)^{s+1}}{(s+1)!}
\end{split}
\end{equation*}
with the last inequality arising from the remainder term in the Taylor
expansion of $e^{\log(C_k\log p_k)}$.  Since from Stirling's approximation $(s+1)! >
((s+1)/e)^{s+1}$ and since from \cite{RosserSchoenfeld} we have
$$P>\frac{1}{e^{\gamma} D_k \log p_k}$$ 
for $k>61$, where
$$D_k = \frac{2 \log^2 p_k}{2 \log^2 p_k - 1}$$ 
 we have
\begin{equation}
\label{PminusT}
P- T^{'}_{s} > \frac{1}{e^{\gamma} D_k \log p_k}- \left(C_k\log p_k
\right) \left( \frac{ e\log(C_k\log p_k) }{s+1}\right)^{s+1}.
\end{equation}

We now wish to find $s$
such that
\begin{equation}
\label{PminusTBound}
\frac{1}{e^{\gamma} D_k \log p_k}- \left(C_k\log p_k \right) \left(
\frac{ e\log(C_k\log p_k) }{s+1}\right)^{s+1} \geq
\frac{1}{2 D_k e^{\gamma} \log p_k}
\end{equation}
Changing variables to 
$$ v = \frac{s+1}{e\log(C_k\log p_k)} $$ 
we can rewrite \eqref{PminusTBound}  as
$$\frac{1}{e^{\gamma} D_k \log p_k}- \frac{C_k\log p_k}{
  \ v^{v\ e\log(C_k\log p_k)}} \geq \frac{1}{2 D_k e^{\gamma}
  \log p_k}$$
and so \eqref{PminusTBound} will hold if
\begin{equation}
\label{inequality2} v^{v\ e\log(C_k\log p_k)} \geq 2 C_k D_k
e^{\gamma}\log^2p_k 
\end{equation}
holds.  
Taking logs on both sides and rearranging \eqref{inequality2} holds if
$$v \log v  \geq \frac{\log\left( 2 C_k D_k e^{\gamma} \log^2 p_k\right)}{e\log(C_k\log p_k)}
$$
holds.
The right-hand side of this inequality is a monotonically decreasing
function of $k$ with the limit $2/e$, and so \eqref{PminusTBound}
holds if we define
$$ v = \frac{ Z\log\left( 2 C_k D_k e^{\gamma} \log^2 p_k\right)}{e\log(C_k\log p_k)}$$
where    $Z=2.159569$  is a solution to 
\begin{equation}
\label{xConstraint}  
Z \log\left(\frac{2 Z}{e}\right) \geq 1
\end{equation}
We thus see that
\eqref{PminusTBound} holds when $s$ is an odd integer such that
\begin{equation}
\label{inequality4}
 s+1 \geq Z\log\left( 2 C_k D_k e^{\gamma} \log^2 p_k\right)
 \end{equation} 
This holds if $s$ is either
$$[Z\log\left( 2 C_k D_k e^{\gamma} \log^2 p_k\right)]$$
or
$$[Z\log\left( 2 C_k D_k e^{\gamma} \log^2 p_k\right)]+1$$
(one of which is odd) and so from \eqref{stevens}, \eqref{PminusT} and \eqref{PminusTBound} we see that 
$$g(n) \leq \left(2  D_k e^{\gamma}  \log p_k \right)k^{Z \log (2 C_k D_k e^{\gamma} )+1}$$
Simplifying and using the fact that 
$$  D_k  \log p_k = k^\frac{\log ( D_k \log p_k)}{\log k}$$
we get
\begin{equation}
\label{fullBound}
g(n) \leq  2 e^{\gamma}  k^{\log (  D_k   \log p_k)/\log k + Z\log(2 C_k D_k e^{\gamma}) +2 Z\log \log p_k + 1}
\end{equation}
Finally we note that all terms of this exponent except $2 Z\log \log p_k$ are either constant or monotonically decreasing with $k$, so if we calculate their values for some fixed $k_0$ and substitute, the bound will hold for all $k \geq k_0$.  We select $k_0=82$  we get
$$g(n) \leq  2 e^{\gamma}  k^{4.78 +2 Z \log \log p_k}$$
for  $k \geq 82$.  Since Kanold's bound is better for $k<82$, this result is valid for all $k > 0$.

Since it is more useful to represent $g(n)$ in terms of $k$ rather than $p_k$, and in terms of integers, we note that for $k > 384$ we have 
$\log p_k < (\log k)^{1.1576}$
and that $2 \times Z \times 1.1576 < 5$, and so we have  
\begin{equation}
\label{simpleBound}
g(n) \leq  2 e^{\gamma}  k^{5 +5 \log \log k}
\end{equation}
for all $k > 384$.  By explicit calculation of the bound in \eqref{fullBound} we find that \eqref{simpleBound} also holds for $k$ in the range $121 \ldots 384$, giving our required result.


\bibliographystyle{amsplain}
 

\end{document}